\documentstyle[11pt,twoside]{article}

\input{mssymb}
\makeatother

\newcommand{\mysec}[2]{%
\section*{\normalsize\hfil\sc {#1}. {#2}\hfill}%
\setcounter{theo}{0}\setcounter{equation}{0}\setcounter{section}{#1}%
\typeout{#1. #2}
\noindent}

\newcommand{\Proof}{\par\noindent{\em Proof. }}
\newcommand{\eop}{\nopagebreak\hspace*{\fill}$\Box$}

\newtheorem{theo}{Theorem}[section]
\newtheorem{lemma}[theo]{Lemma}
\newtheorem{cor}[theo]{Corollary}

\newtheorem{prop}[theo]{Proposition}
\newtheorem{definition}[theo]{Definition}

\newcounter{abc}   % Counter fr statements-environment wird deklariert
\newcounter{iiiii} % Counter fr aequivalenz-environment wird deklariert

\newenvironment{aequivalenz}
{\setcounter{iiiii}{0}
\begin{list}%
{{\rm (\roman{iiiii})}}%  Falls die items nicht angegeben sind: i)u.s.w.
{\usecounter{iiiii}
\parsep=0pt plus 1pt
\topsep=1pt plus 2pt minus 1pt
\itemsep=1pt plus 2pt minus 1pt
\leftmargin=3\baselineskip
\labelsep=.6\baselineskip
\labelwidth=2.4\baselineskip
\rightmargin 0pt}%
}%               Das war das zweite Argument von "newenvironment"
{\end{list}}

\newenvironment{statements}%
{\setcounter{abc}{0}
\begin{list}%
{{\rm (\alph{abc})}}%  Falls die items nicht angegeben sind: (a) u.s.w.
{\usecounter{abc}
\parsep=0pt plus 1pt
\topsep=1pt plus 2pt minus 1pt
\itemsep=1pt plus 2pt minus 1pt
\leftmargin=3\baselineskip
\labelsep=.6\baselineskip
\labelwidth=2.4\baselineskip
\rightmargin 0pt}%
}%               Das war das zweite Argument von "newenvironment"
{\end{list}}

{\begin{list}%
{$\bullet$}%  Falls die items nicht angegeben sind: Punkte
{\leftmargin=5\baselineskip
\labelsep=\baselineskip
\labelwidth=2.5\baselineskip
\parsep=0pt plus 1pt
\topsep=1pt plus 2pt minus 1pt
\itemsep=1pt plus 2pt minus 1pt
\rightmargin 0pt}%
}%               Das war das zweite Argument von "newenvironment"
{\end{list}}

\newif\ifrefsc
\let\thebibliographyalt=\thebibliography                                %
\def\thebibliography#1                                                  %
 {\ifrefsc
 \section*{\normalsize\hfil\sc References\hfill}
 \small
 \list                    %
 {[\arabic{enumi}]}{\settowidth\labelwidth{[#1]}\leftmargin\labelwidth  %
 \advance\leftmargin\labelsep                                           %
 \usecounter{enumi}}                                                    %
 \def\newblock{\hskip .11em plus .33em minus .07em}                     %
 \sloppy\clubpenalty4000\widowpenalty4000\vfuzz 2pt                      %
 \sfcode`\.=1000\relax                                                  %
 \else\thebibliographyalt{#1}\fi}                                         %

\makeatletter
\def\ps@smallheadings{            % Nur bei doppelseitigem Ausdruck
  \def\@oddfoot{}                 % Aufruf mit: pagestyle{smallheadings}
  \def\@evenfoot{}                % keine Footer
  \def\@evenhead{\hbox to \textwidth {%
  \vbox{\hbox to \textwidth 
        {\thepage\hfil\strut{\footnotesize\sc\Autor}\hfil}\vss}}}
  \def\@oddhead{\hbox to \textwidth {%
\vbox{\hbox to \textwidth 
        {\hfil{\footnotesize\sc\Kurztitel}\hfil\strut \thepage}\vss}}}
    }
\makeatother

\makeatletter
\def\eqalignno#1{\displ@y \tabskip\@centering
  \halign to\displaywidth{\hfil$\@lign\displaystyle{##}~$\tabskip\z@skip
    &$\@lign\displaystyle{{}##}$\hfil\tabskip\@centering
    &\llap{$\@lign##$}\tabskip\z@skip\crcr
    #1\crcr}}

\makeatother

\def\ersteSeite{\vspace*{32pt plus 2pt minus 2pt}\begin{center}
{\Large\sf\Titel}\\[15pt]{\sc\Autor}\\[26pt plus 2pt minus 2pt]
\end{center}\Abstrakt\vspace{0pt plus 2pt }\thispagestyle{empty}}

\def\Abstrakt{\begin{quote}\small\noindent{\sc Abstract.}
\Abstrakttext\end{quote}}

\newcommand{\N}{{\Bbb N}}
\newcommand{\R}{{\Bbb R}}

\newcommand{\T}{{\Bbb T}}
\newcommand{\Z}{{\Bbb Z}}

\newcommand{\Id}{I\mkern-1mud}

\newcommand{\alp}{\alpha}
\newcommand{\gam}{\gamma}
\newcommand{\Gam}{\Gamma}
\newcommand{\del}{\delta}
\newcommand{\eps}{\varepsilon}
\newcommand{\lam}{\lambda}
\newcommand{\Lam}{\Lambda}

\newcommand{\Ome}{\Omega}

\newcommand{\klam}{\bigl(}
\newcommand{\mer}{\bigr)}

\newcommand{\loglike}[1]{\mathop{\rm #1}\nolimits}

\newcommand{\ex}{\loglike{ex}}

\newcommand{\lin}{\loglike{lin}}  

\newcommand{\re}{\loglike{Re}}

\newcommand{\ran}{\loglike{ran}}

\newcommand{\qqfa}{\qquad\forall}

\newcommand{\bea}{\begin{eqnarray*}}
\newcommand{\eea}{\end{eqnarray*}}
\newcommand{\beq}{\begin{equation}}
\newcommand{\eeq}{\end{equation}}
\newcommand{\begsta}{\begin{statements}}
\def\endsta{\end{statements}}
\newcommand{\begaeq}{\begin{aequivalenz}}
\def\endaeq{\end{aequivalenz}}

\newcommand{\iy}{\infty}
\newcommand{\dopu}{{:}\allowbreak\ }

\newcommand{\kle}{{\ell^{1}}}

\newcommand{\rest}[2]{#1\raisebox{-0.3ex}{\mbox{$\mid_{#2}$}}} 

\newcommand{\DE}{Daugavet equation}
\def\Remark{2.5}%
\def\Example{3.2}

\refsctrue

\begin{document}           \tracingpages=1

%\preliminaryversion

\def\Titel{The Daugavet equation for operators on function spaces}
\def\Kurztitel{\Titel     }
\def\Autor{Dirk Werner}
\def\Abstrakttext{
We prove the norm identity $\|\Id+T\| =1+\|T\|$, which is known as
the Daugavet equation, for weakly compact operators $T$ on natural
function spaces such as function algebras and $L^{1}$-predual spaces,
provided a non-discreteness assumption is met. We also consider
$c_{0}$-factorable operators and operators on
$C_{\Lambda}$-spaces.
}

\ersteSeite

\mysec{1}{Introduction}%

In his 1963 paper \cite{Daug} Daugavet proved the remarkable norm identity
\beq\label{eqx.y}
\|\Id + T\| = 1 + \|T\|
\eeq
for a compact operator on $C[0,1]$; in the sequel (\ref{eqx.y}) has
become known as the {\em \DE}. (\ref{eqx.y}) has proved useful in
approximation theory; Ste\v{c}kin used it in order to provide
precise lower bounds for trigonometric approximations \cite{Steckin}.
Daugavet's result was extended to other classes of operators on
spaces of continuous functions $C(S)$, notably to weakly compact
operators \cite{FoiSin} and also to operators on $L^{1}$-spaces; we
refer to \cite{Abra2} and 
\cite{AbraAB} for a more detailed account of the history of
the subject. 

In this paper we suggest a systematic, yet simple approach to
studying the \DE\ for operators on subspaces of $C(S)$-spaces. Since
every Banach space embeds isometrically into a $C(S)$-space, some
restriction on the embedding has to be imposed. Also, a natural
obstruction for the \DE\ to hold for  -- say -- compact operators on
$C(S)$ is the presence of isolated points in $S$, because an isolated
point of $S$ immediately gives rise to a one-dimensional operator
(i.e., an operator with one-dimensional range) $T$ on $C(S)$ with
$\|\Id -T\|=1$. In section~2 we formulate conditions on an isometric
embedding $J$ of a Banach space $X$ into $C(S)$, labelled (N1) and
(N2), which we will assume in order that $X$ be ``nicely embedded;''
in addition we will require a non-discreteness condition~(N3).

In the next section we will present a necessary and sufficient
condition for the  \DE\ on a nicely embedded space, and we check it
for weakly compact operators and for operators factoring through a
subspace of $c_{0}$; note that both these classes encompass the
compact operators. 
Whereas section~2 has a deliberately technical flavour, we give
applications to concrete spaces in section~3. In particular, we deal
with function algebras, $L^{1}$-predual spaces and translation
invariant spaces, and we establish the \DE\ for various classes of
operators.

Condition (N2) will be given in terms of the $L$-structure of
$X^{*}$. We recall the relevant definitions. A closed subspace $F$ of a
Banach space $E$ is an {\em $L$-summand}\/ if there is a projection
$\Pi$ from $E$ onto $F$ such that
$$
\|\xi\| = \|\Pi(\xi)\| + \|\xi - \Pi(\xi)\| \qqfa \xi\in E.
$$
Dual to this notion is the definition of an $M$-ideal: $F\subset E$
is an {\em $M$-ideal}\/ if its annihilator $F^{\bot}\subset E^{*}$ is
an $L$-summand. These concepts are studied in detail in \cite{HWW}.
Roughly speaking, our conditions (N1)--(N3) mean that the function
space $X$ has a rich and nondiscrete $M$-ideal structure.

We use standard notation such as $L(X)$ for the space of all bounded
linear operators on a Banach space $X$, $B_{X}$ for 
the closed unit ball of $X$ and
$\ex C$ for the set of extreme points of a convex set~$C$.

\mysec{2}{Nicely embedded Banach spaces}%

In order to formulate the lemmas of this section succinctly, we need
to introduce some vocabulary. Let $S$ be a Hausdorff
topological space, and let
$C^{b}(S)$ be the sup-normed Banach space of all bounded continuous
scalar-valued functions. The functional $f\mapsto f(s)$ on $C^{b}(S)$
is denoted by $\del_{s}$. We say that a linear map $J\dopu X \to
C^{b}(S)$ on a Banach space $X$ is a {\em nice embedding\/} and that
$X$ is {\em nicely embedded\/} into $C^{b}(S)$ if $J$ is an isometry
such that for all $s\in S$ the following properties are satisfied:
\smallskip
\begsta\em
\item[\rm(N1)]
For $p_{s}:= J^{*}(\del_{s})\in X^{*}$ we have $\|p_{s}\|=1$.
\item[\rm(N2)]
$\lin\{p_{s}\}$ is an $L$-summand in $X^{*}$.
\endsta
\smallskip
The latter condition can also be rephrased by saying that the kernel
of $p_{s}$ is an $M$-ideal.    
We will discuss examples of nicely embedded Banach spaces in
section~3.

Throughout this section we will stick to the following notation. Let
$J\dopu X\to C^{b}(S)$ be a nice embedding, $p_{s}=J^{*}(\del_{s})$,
and let $T\in L(X)$.
We put
$$
q_{s}:= (JT)^{*}(\del_{s}) = T^{*}(p_{s})\in X^{*}
$$
and note that $s\mapsto q_{s}$ is weak$^{*}$ continuous and $\|T\|=
\sup_{s}\|q_{s}\|$. Likewise, $s\mapsto p_{s}$ is weak$^{*}$
continuous. By~(N2) there is a family of projections $\Pi_{s}$, $s\in
S$, with $\ran \Pi_{s}=\lin\{p_{s}\}$ such that
$$
\|x^{*}\| = \|\Pi_{s}(x^{*})\| + \|x^{*}-\Pi_{s}(x^{*})\| \qqfa x^{*}\in
X^{*}
$$
and a family of functionals $\pi_{s}\in X^{**}$, $s\in S$, such that
$$
\Pi_{s}(x^{*}) = \pi_s(x^{*}) p_{s} \qqfa x^{*}\in X^{*}.
$$
In particular, we have $\pi_{s}(p_{s})=1$.

We will also need the equivalence relation
\beq\label{eq2.2}
s \sim t \ \mbox{ if and only if }\  \Pi_{s} =\Pi_{t}
\eeq
on $S$. Then $s$ and $t$ are equivalent if and only if $p_{s}$ and
$p_{t}$  are linearly dependent, which implies by (N1) that
$p_{t}=\lam p_{s}$ for some scalar of modulus~1. The equivalence
classes of this relation are obviously closed.

In some of the lemmas to follow we will additionally have to assume
\smallskip
\begsta\em
\item[\rm (N3)]
None of the equivalence classes $Q_{s}=\{t\in S\dopu s\sim
t\}$ contains an interior point.
\endsta
\smallskip
If the set $\{p_{s}\dopu s\in S\}$ is linearly independent, this
simply means that
\smallskip
\begsta\em
\item[\rm (N$3'$)]
$S$ does not contain an isolated point.
\endsta
\smallskip
By (N2), the $p_{s}$ are linearly independent as soon as they are
pairwise linearly independent.

We now give the basic criterion for an operator on a nicely embedded
Banach space to satisfy the \DE.

\begin{lemma}\label{2.1}
Let $X$ be nicely embedded into $C^{b}(S)$ so that\/ {\rm (N1)} and\/
{\rm (N2)} are valid, and let $T\in L(X)$. For $\eps>0$ put 
$U_{\eps}= \{s\in S\dopu \|q_{s}\|>\|T\|-\eps\}$, which is an open
subset of $S$. Then $T$ satisfies the \DE\/~{\rm(\ref{eqx.y})} if and
only if
\beq\label{eq2.1}
\sup_{s\in U_{\eps}}  \klam |1+\pi_{s}(q_{s})| - (1+|\pi_{s}(q_{s})|)
\mer \ge0 \qqfa \eps>0.
\eeq
\end{lemma}

\Proof
First assume (\ref{eq2.1}).
We observe
\bea
\|p_{s}+q_{s}\| &=& \|\Pi_{s}(p_{s}+q_{s})\|
                 + \| (p_{s}+q_{s}) - \Pi_{s}(p_{s}+q_{s})\| \\
&=& |1+\pi_{s}(q_{s})| + \|(\Id-\Pi_{s})(q_{s})\|
\eea
and
$$
1+\|q_{s}\| = 1+ |\pi_{s}(q_{s})| + \|(\Id-\Pi_{s})(q_{s})\|
$$
for all $s\in S$. Applying the assumption, for some $\eps>0$, 
we obtain from this
\bea
\|\Id+T\| &=&
\sup_{s\in S} \|p_{s}+q_{s}\|
~\ge~ \sup_{s\in U_{\eps}} \|p_{s}+q_{s}\| \\
&=&
\sup_{s\in U_{\eps}} \klam 
|1+\pi_{s}(q_{s})| + \|(\Id-\Pi_{s})(q_{s})\| \mer  \\
&\ge&
1+ |\pi_{s}(q_{s})| -\eps + \|(\Id-\Pi_{s})(q_{s})\| \\
&&\qquad \mbox{for some }s\in U_{\eps} \\
&=&
1 + \|q_{s}\| -\eps \\
&>&
1+\|T\| -2 \eps \qquad(\mbox{since }s\in U_{\eps}).
\eea
This proves (\ref{eqx.y}), since $\eps>0$ was arbitrary.

The proof of the converse implication follows the same lines.
\eop

\begin{cor}\label{2.1a}
If $X$ is nicely embedded into $C^{b}(S)$ so that\/ {\rm (N1)} and\/
{\rm (N2)} are valid and if $T\in L(X)$, then there is a scalar
$\lam$, $|\lam|=1$, such that $\lam T$
satisfies the \DE\/~{\rm(\ref{eqx.y})}.
\end{cor}

In the case of real Banach spaces (\ref{eq2.1}) is equivalent to
$$
\sup_{s\in U_{\eps}} \pi_{s}(q_{s}) \ge0 \qqfa \eps>0,
$$
and Corollary~\ref{2.1a} simply says that $T$ or $-T$ satisfies the
\DE.

It remains to give examples of classes of operators for which
(\ref{eq2.1}) is valid. This will be done in Lemmas~\ref{2.3}
and~\ref{2.5}. But first we single out a simple estimate that will be
used in the proof of those lemmas.

\begin{lemma}\label{2.2}
Suppose $X$ is nicely embedded into $C^{b}(S)$ such that\/ {\rm(N1)}
and\/ {\rm(N2)} hold. If $t_{1},\ldots,t_{k}$  are pairwise
nonequivalent points $($for the equivalence relation $\sim$ of\/
{\rm(\ref{eq2.2}))}, then
$$
\|x^{*}\| \ge \sum_{j=1}^{k} \|\Pi_{t_{j}}(x^{*})\| \qqfa x^{*}\in
X^{*}.
$$
\end{lemma}

\Proof
Let $\Pi=\sum_{j=1}^{k}\Pi_{t_{j}}$ be the $L$-projection with range
$\lin\{p_{t_{1}},\ldots,p_{t_{k}}\}$. Then
$$
\|x^{*}\| \ge \|\Pi(x^{*})\| = 
\sum_{j=1}^{k} \|\Pi_{t_{j}}(x^{*})\|,
$$
which implies our claim.
\eop

\begin{lemma}\label{2.3}
Suppose $X$ is nicely embedded into $C^{b}(S)$ such that\/ {\rm(N1),
(N2)}
and\/ {\rm(N3)} hold, and let $T\in L(X)$ be an operator. If
\beq\label{eq2.3}
s\mapsto \pi_{t}(q_{s}) \mbox{ is continuous for all }t\in S,
\eeq 
then $T$
satisfies\/ {\rm(\ref{eq2.1})} and consequently the \DE\/~{\rm(\ref{eqx.y})}.
\end{lemma}

\Proof
We consider the function $f$ on $S\times S$ defined by
$f(s,t)=\pi_{t}(q_{s})$. Then our assumption on $T$ means that
$s\mapsto f(s,t)$ is continuous for all $t\in S$. Now we check
condition (\ref{eq2.1}) and argue by contradiction. If (\ref{eq2.1})
were false, we would find an open set $U\neq\emptyset$ and some
$\beta>0$  such that
$$
|1 + f(s,s)| - (1+|f(s,s)|) <-2\beta \qqfa s\in U.
$$
In particular,
\beq\label{eq2.4}
|f(s,s)| > \beta \qqfa s\in U.
\eeq

Let $s_{1}\in U$  be arbitrary. By continuity of $f $ in the first
variable, there is an open neighbourhood $U_{1}\subset U$ of $s_{1}$ 
such that
\beq\label{eq2.5}
|f(u,s_{1})| >\beta \qqfa u\in U_{1}.
\eeq
Since the equivalence class $Q_{s_{1}}$ does not contain interior
points, we may pick some $s_{2}\in U_{1}$, $s_{2}\notin Q_{s_{1}}$.
Consequently, we have by (\ref{eq2.4}) and (\ref{eq2.5}) 
$$
|f(s_{2},s_{2})| > \beta,\quad |f(s_{2},s_{1})| > \beta.
$$
We proceed to find an open neighbourhood $U_{2}\subset U_{1}$ of $s_{2}$ 
such that
\beq\label{eq2.6}
|f(u,s_{2})| >\beta \qqfa u\in U_{2}
\eeq
and some $s_{3}\in U_{2}$, $s_{3}\notin Q_{s_{1}}\cup Q_{s_{2}}$,
with
$$
|f(s_{3},s_{3})| > \beta,\quad |f(s_{3},s_{2})| > \beta,
\quad |f(s_{3},s_{1})| > \beta
$$
by (\ref{eq2.4}), (\ref{eq2.6}) and (\ref{eq2.5}).

Continuing in the obvious manner, we arrive at a sequence
$s_{1},s_{2},\ldots$ in $S$ such that for each $k\in \N$
$$
|f(s_{k},s_{j})| > \beta \qqfa j=1,\ldots,k.
$$
But by Lemma~\ref{2.2}, since the $\Pi_{s_{1}}, \Pi_{s_{2}}, \ldots$
are different $L$-projections, we have for each $k\in \N$
$$
\|T\| \ge \|q_{s_k}\| \ge \sum_{j=1}^{k} \|\Pi_{s_{j}}(q_{s_{k}})\|
= \sum_{j=1}^{k} |f(s_{k},s_{j})| > k\beta,
$$
which clearly contradicts the continuity of $T$.

Thus, the lemma is proved.
\eop
\bigskip
\addtocounter{theo}{1}

\noindent
{\em Remark \thetheo\/} \ \label{Remark}
Weakly compact operators fulfill the continuity assumption
(\ref{eq2.3}) in
Lemma~\ref{2.3}. In fact, if $T\in L(X)$ is weakly compact, then
$T^{*}$ is weak$^{*}$-weakly-continuous, and $s\mapsto
q_{s}=T^{*}p_{s}$ is weakly continuous since $s\mapsto p_{s}$ is
weak$^{*}$ continuous. We will meet non-weakly-compact operators that
fulfill (\ref{eq2.3}) in Example~\Example.
\bigskip

We will now discuss a different class of operators
satisfying~(\ref{eq2.1}). We say that $T\in L(X)$ factors through a
subspace of $c_{0}$ if there is a closed subspace $E\subset c_{0}$
together with continuous operators $T_{1}\dopu X\to E$  and
$T_{2}\dopu E\to X$ such that $T=T_{2}T_{1}$. As already noted in the
introduction, 
every compact operator has this property.

\begin{lemma}\label{2.5}
Suppose $X$ is nicely embedded into $C^{b}(S)$ such that\/ {\rm(N1),
(N2)}
and\/ {\rm(N3)} hold, and let $T\in L(X)$ be an operator
factoring through a subspace of $c_{0}$. Assume further that $S$ is a
Baire space. Then $T$
satisfies\/ {\rm(\ref{eq2.1})} and consequently the \DE\/~{\rm(\ref{eqx.y})}.
\end{lemma}

\Proof
As indicated above, let us write
$$
JT\dopu X \stackrel{T_{1}}{\longrightarrow} E
\stackrel{T_{2}}{\longrightarrow} X
\stackrel{J}{\longrightarrow} C^{b}(S).
$$
We define functionals $x_{n}^{*}\in X^{*}$ $(n\in \N)$ and
$a_{s}^{*}\in E^{*}$ $(s\in S)$ by
$$
\begin{array}{rcll}
x_{n}^{*}(x) &=& (T_{1}x)(n) \quad&(x\in X),\\[2pt]
a_{s}^{*}(a) &=& (JT_{2}a)(s)     &(a\in E).
\end{array}
$$
Observe that $\sup_{n}\|x_{n}^{*}\| = \|T_{1}\|$ and $\sup_{s}
\|a_{s}^{*}\| = \|T_{2}\|$. Let $\nu_{s}\in \kle$ be a Hahn-Banach
extension of $a_{s}^{*}\in E^{*}$. Then we have
$$
(JTx)(s) = a_{s}^{*}(T_{1}x) = \sum_{n=1}^{\iy}
\nu_{s}(n)x_{n}^{*}(x) = \biggl( \sum_{n=1}^{\iy } \nu_{s}(n)x_{n}^{*}
\biggr)(x)
$$
and consequently
$$
q_{s}= \sum_{n=1}^{\iy } \nu_{s}(n)x_{n}^{*};
$$
note that this series is absolutely norm-convergent.

In order to achieve the proof proper of Lemma~\ref{2.5}, we define
$$
S' = \{t\in S\dopu \pi_{t}(q_{s}) =0 \ \forall s\in S\},
$$
and we claim that $S'$ is dense in $S$, which clearly
implies~(\ref{eq2.1}). In fact, from
$$
\Pi_{t}(q_{s}) = \sum_{n=1}^{\iy} \nu_{s}(n) \Pi_{t}(x_{n}^{*})
$$
we infer that
$$
S'':= \{ t\in S\dopu \Pi_{t}(x_{n}^{*})=0 \ \forall n\in\N\} \subset
S',
$$
and
$$
S{\setminus}S'' = \bigcup_{n\in\N}  \{ t\in S\dopu \Pi_{t}(x_{n}^{*})
\neq0 \}.
$$
But by Lemma~\ref{2.2}, each set $\{ t\dopu
\Pi_{t}(x_{n}^{*})\neq0\}$ consists of at most countably many
equivalence classes for $\sim$, for there are at most $\|T_{1}\|/\del$ many
nonequivalent $t$ with $\|\Pi_{t}(x_{n}^{*})\|\ge\del$.
Therefore $S{\setminus}S''$ is a countable union of (by
assumption~(N3)) nowhere dense sets. The Baire property yields that
$S''$ and hence $S'$ is dense.
\eop

\mysec{3}{Applications}%

Now we will apply the results of the previous section to some natural
classes of function spaces. The most obvious example of a nicely
embedded space is of course $C_{0}(S)$, $S$ a locally compact
Hausdorff space, with $J=$ the natural inclusion into $C^b(S)$.

\begin{prop}\label{3.1}
Let $S$ be a locally compact Hausdorff space without isolated points.
If $T\in L(C_{0}(S))$ is weakly compact or factors through a subspace
of $c_{0}$, then $T$ satisfies the \DE\/~{\rm(\ref{eqx.y})}. More
generally, it is enough that the function $s\mapsto
(T^{*}\del_{s})(\{t\})$ on $S$ is continuous for all $t\in S$  
in order that $T$ satisfies the \DE.
\end{prop}

\Proof
Since here $s\sim t$ iff $s=t$, we see that (N3$'$) and hence (N3)
are fulfilled so that the assertion follows from Lemma~\ref{2.3},
Remark~\Remark\  and Lemma~\ref{2.5}. (Observe that $\Pi_{t}$ is the
$L$-projection $\mu\mapsto \mu(\{t\})\del_{t}$ on $M(S)\cong
(C_{0}(S))^{*}$ in the present context so that $\pi_{t}(q_{s}) =
(T^{*}\del_{s})(\{t\})$.)
\eop
\bigskip

Most of Proposition~\ref{3.1} is already known; we refer to
\cite{Dirk8} for
a straightforward proof and to \cite{AbraAB}, \cite{Ansari},
\cite{FoiSin} and \cite{Holub2} for different other approaches.
\addtocounter{theo}{1}

\bigskip
\noindent
{\em Example \thetheo\/} \
\label{Example}
There are non-weakly-compact operators on $C(\T)$ such that $s\mapsto
(T^{*}\del_{s})(\{t\})=0$ for all $t\in\T$; a fortiori these are
continuous functions. In fact, every convolution operator $T\dopu
f\mapsto f \ast \mu$ for a continuous ($=$ diffuse) singular measure
has this property: In this case $(T^{*}\del_{s})(\{t\})=\mu(\{
st^{-1}\})=0$, since $\mu $ is continuous; 
and a result due to Cost\'e \cite[p.~90]{DiUh} implies that $\mu $
would be absolutely continuous if $T$ were weakly compact. I am
grateful to W.~Hensgen for pointing out this example to me.

Therefore we see that such convolution operators satisfy the \DE, a
fact that can also be checked directly.

\bigskip
We now turn to algebras of functions. A {\em function algebra\/} $A$
on a compact Hausdorff space $K$ is a closed subalgebra of the space
of complex-valued functions $C(K)$ separating the points of $K$ and
containing the constant functions. To each function algebra $A$ on $K$
one can associate a distinguished subset $\partial A\subset K$,
called the {\em Choquet boundary}, defined by
$$
\partial A = \{k\in K\dopu \rest{\del_k}{A} 
\mbox{ is an extreme point of }B_{A^{*}} \}.
$$
This notion will be instrumental in the proof of the following
result.

\begin{theo}\label{3.3}
Let $A$ be a function algebra such that its Choquet boundary
$\partial A$ does not contain an isolated point. Then every operator
$T\in L(A)$ which is weakly compact or factors through a subspace
of $c_{0}$ satisfies the \DE\/~{\rm(\ref{eqx.y})}.
\end{theo}

\Proof
We will verify  that $A$ is nicely embedded into $C^{b}(\partial A)$
such that (N1), (N2) and (N3) are valid. We consider the mapping
$J\dopu A\to C^{b}(\partial A)$, $Jf= \rest{f}{\partial A}$. It is a
well-known consequence of the Krein-Milman theorem that $J$
is an isometry, cf.\ \cite[p.~180f.]{Choq2}\ for details. 
Since $J^{*}(\del_{k})= \rest{\del_k}{A}$,
condition (N1) is fulfilled by construction, and (N2) is a result due
to Hirsberg \cite{Hir}, see also \cite[p.~15 and Th.~V.4.2]{HWW}. As in
Proposition~\ref{3.1}, (N3) reduces to (N3$'$) which is part of the
assumption of Theorem~\ref{3.3}.
We finally observe that $\partial A$ is homeomorphic to the extreme
boundary of $\{\ell\in A^{*}\dopu \|\ell\| = |\ell({\bf1})| = 1\}$
(the state space of $A$). By a theorem of Choquet's
\cite[p.~146]{Choq2} $\partial A$
is a Baire space.

Hence we can apply Lemma~\ref{2.3}, Remark~\Remark\ and
Lemma~\ref{2.5} to finish the proof of Theorem~\ref{3.3}.
\eop

\begin{cor}\label{3.3a}
Let $A$ be a function algebra such that its Choquet boundary
$\partial A$ does not contain an isolated point. Suppose $\{0\}\neq
E\subset A$ is either reflexive or isomorphic to $c_{0}$, and suppose
that $E$ is the kernel of a projection~$P$. Then $\|P\|\ge2$.
\end{cor}

\Proof
The operator $\Id-P$ is either weakly compact or factors through
$c_{0}$, and it is a nonzero projection. Thus, by Theorem~\ref{3.3},
$$
\|P\| = \|\Id - (\Id-P)\| = 1 + \|\Id-P\| \ge2. \eqno\Box
$$

\bigskip
In particular, finite-codimensional proper subspaces are complemented
only by projections of norm~$\ge2$.

A similar corollary can be formulated for other classes of Banach
spaces discussed in this section; but note that the only complemented
reflexive subspaces of $C(K)$ (or of an $L^{1}$-predual space, see
below) are finite-dimensional, by the Dunford-Pettis property of
those spaces.

The \DE\ for weakly compact operators on function algebras was first
established by Wojtaszczyk \cite{Woj92} using a different argument.
He also remarks that $\partial A$ fails to contain an isolated point
if $A$ fails to contain nontrivial idempotents -- a property shared
by the disk algebra and the algebra of bounded analytic functions on
the unit disk.

Again, it would have been enough to require (\ref{eq2.3}) instead of
the weak compactness of $T$. 
%We give a natural example where this
%generalisation is useful.
%\addtocounter{theo}{1}
%
%\bigskip
%\noindent
%{\em Example~\thetheo}\/ \ 
%\label{Exampletwo}
%Let $A$ be the disk algebra, considered as a subspace of $C(\T)$. If
%$\mu$ is a continuous singular measure on $\T$, then the convolution
%operator $T\dopu f\mapsto f \ast\mu$  maps $A$ into $A$, satisfies
%the \DE, but is not weakly compact. In fact, we have the
%decomposition $A^{*}= M_{\rm sing} \oplus L^{1}(\T)/H_{0}^{1}$ by the
%F.~and M.~Riesz theorem, and for the convolution operator $S\dopu
%f\mapsto f\ast\mu$ on $C(\T)$ we have $(S^{*}\del_{s})(B)= \mu(s-B)$
%(we preferred to write the group operation on $\T$ as addition) for
%every Borel set $B$. Since $\mu$ is singular, we can identify
%$T^{*}(\rest{\del_{s}}{A})$ with this measure, for each $s\in \T$.
%Thus $T^{*}(\rest{\del_{s}}{A})\in M_{\rm sing}$ in a canonical way,
%and $\pi_{t}(q_{s}) = (T^{*}\rest{\del_{s}}{A})(\{t\}) = 0$ for all
%$t$ so that (\ref{eq2.3})  is fulfilled and $T$ satisfies the \DE.
%In Example~\Example\ we noticed that $S$ is not weakly compact, hence
%$s\mapsto (S^{*}\del_{s})(B)$ is not continuous for some Borel set $B$
%\cite[p.~153f.]{DiUh}. Consequently, $s\mapsto q_{s} = 
%T^{*}(\rest{\del_{s}}{A})$ is not weakly continuous either, and $T$
%is not weakly compact.
%
%We will resume the discussion of this example in
%Example~\Examplethree.
%\bigskip

The next class of Banach spaces we wish to investigate are the
$L^{1}${\em -predual spaces}\/ $X$ defined by the requirement that
$X^{*}$ is isometric to a space of integrable functions $L^{1}(\mu)$.
We consider the equivalence relation $p \sim q$ iff $p$ and $q$ are
linearly dependent  on  $\ex B_{X^{*}}$, and we equip the quotient
space $\ex B_{X^{*}}/{\sim}$ with the quotient topology of the
weak$^{*}$  topology. Of course, the following result contains
Proposition~\ref{3.1} as a special case.

\begin{theo}\label{3.4}
Let $X$ be an $L^{1}$-predual space such that $\ex B_{X^{*}}/{\sim}$
does not contain an isolated point. Then every operator
$T\in L(X)$ which is weakly compact or factors through a subspace
of $c_{0}$ satisfies the \DE\/~{\rm(\ref{eqx.y})}.
\end{theo}

\Proof
Let $J\dopu X\to C^{b}(\ex B_{X^{*}})$ be the canonical isometry. For
$s\in \ex B_{X^{*}}$ we have $J^{*}(\del_{s})=s$; hence (N1) and (N2)
are satisfied. (Observe that the linear span of an extreme point in
an $L^{1}$-space is an $L$-summand.) Also, (N3) holds by assumption
on $\ex B_{X^{*}}$. 
Finally, we again invoke Choquet's theorem \cite[p.~146]{Choq2} 
to ensure that $\ex
B_{X^{*}}$ is a Baire space. Thus, Lemma~\ref{2.3}, Remark~\Remark\ 
and Lemma~\ref{2.5} yield Theorem~\ref{3.4}.
\eop

\bigskip
\addtocounter{theo}{1}
\noindent
{\em Example~\thetheo}\/ \ 
\label{Examplefour}
Let $\Ome\subset \R^{d}$  be open and bounded and consider the
sup-normed space 
$$
H(\Ome)= \{f\in C(\overline{\Ome})\dopu f \mbox{ is harmonic on }
\Ome\}.
$$
This is an order unit space whose state space $K$ is a Choquet
simplex so that $H(\Ome)$ is an $L^{1}$-predual space (cf.\
\cite{EfKaz}). The extreme points of $K$, that is up to
scalar multiples the extreme points of 
$B_{H(\Ome)^{*}}$, can be identified with the set $\partial_{r}\Ome$
of regular boundary points (in the sense of potential theory) of~$\Ome$. 
Thus we conclude that every operator on $H(\Ome)$ which is
weakly compact or factors through a subspace of $c_{0}$  satisfies
the \DE, provided $\partial_{r}\Ome$ does not contain isolated
points. It is classical that  $\partial_{r}\Ome= \partial\Ome$ if the
boundary of $\Ome$ is sufficiently smooth or if $\Ome$ is simply
connected and ${d=2}$.

\bigskip
Our final application deals with translation invariant spaces. Let
$G$ be an infinite compact abelian group with dual group
$\Gamma$ and Haar measure~$m$; we denote the group operation in $\Gam$
by~$+$. For $\Lam\subset \Gam$ the space of
$\Lam$-spectral continuous functions is defined by
$$
C_{\Lam} = \{f\in C(G)\dopu \widehat{f}(\gam)=0 \ \forall
\gam\notin\Lam\},
$$
where $\widehat{\gam}$ denotes the Fourier transform of $f$. These
spaces are known to be precisely the closed translation invariant subspaces
of $C(G)$. Likewise, one defines spaces of $\Lam$-spectral measures
$M_{\Lam}$ and $\Lam$-spectral integrable functions $L^{1}_{\Lam}$.

A subset $\Lam\subset \Gam$ is called a {\em Riesz set}\/ if
$M_{\Lam}\subset L^{1}(m)$; the chief example of a Riesz subset of
$\widehat{\T}=\Z$ is $\N$. For an in-depth analysis of this class
of sets and its relation to Banach space geometry we refer to
\cite{G-Riesz}, see also Chapter~IV.4 in \cite{HWW}. Here we consider
a broader class of sets which we propose to call semi-Riesz sets: If
$M_{\rm diff}$ denotes the space of diffuse ($=$ continuous)
measures on $G$,
i.e., those which map singletons to $0$, then $\Lam\subset \Gam$ is a
{\em semi-Riesz set}\/ if $M_{\Lam}\subset M_{\rm diff}$. Obviously,
Riesz sets are semi-Riesz,
but there are others; typical examples of proper semi-Riesz
sets are  spectra of Riesz products. To be definite, consider the
Riesz product
$\mu= w^{*}\mbox{-}\lim_{n\to\infty} \prod_{k=0}^{n} (1+\cos
4^{k}t)\,dm \in M[0,2\pi)\cong M(\T)$.
Let $\Lambda= \{ \sum_{k=0}^{n} \eps_{k} 4^{k}\dopu \eps_{k}= -1,0,1,\ 
n\in\N\}$; then $\mu\in M_{\Lambda}(\T)$, and $\mu$ is not absolutely
continuous. So $\Lambda $ is not a Riesz set. It is, however,
semi-Riesz, which can be deduced from a theorem of Wiener's (cf.\
\cite[p.~415]{GraMcG}). In fact, for $\nu\in M_{\Lambda}(\T)$ we have
$$
\frac{1}{2N+1} \sum_{k=-N}^{N} |\widehat{\nu}_{k}|^{2} \le
\|\widehat{\nu}\|_{\infty}^{2} \frac{\# \{\lambda\in\Lambda\dopu
|\lambda|\le N\}}{2N+1} \to0.
$$
Wiener's theorem implies that $\nu$ is diffuse.

\begin{theo}\label{3.5}
Let $G$ be a compact abelian group and suppose $\Lam$ is a subset of
the dual group $\Gam$ such set $\Gam{\setminus}(-\Lam)$ is a
semi-Riesz set. Then every operator
$T\in L(C_{\Lam})$ which is weakly compact\/ $($or merely satisfies\/
{\rm (\ref{eq2.3}))} or factors through a subspace
of $c_{0}$ satisfies the \DE\/~{\rm(\ref{eqx.y})}.
\end{theo}

\Proof
Let $J\dopu C_{\Lam}\to C(G)$ be the identical mapping. Then $J^{*}$
is the quotient map onto
\bea
C_{\Lam}^{*} ~\cong~ M(G)/(C_{\Lam})^{\bot}
&\cong& M(G)/M_{\Gam{\setminus}(-\Lam)} \\
&\cong& (\kle(G) \oplus_{1} M_{\rm diff})/M_{\Gam{\setminus}(-\Lam)} \\
&\cong&
\kle(G) \oplus_{1} M_{\rm diff}/M_{\Gam{\setminus}(-\Lam)} 
\eea
with $\oplus_{1}$ denoting $\kle$-direct sums, by assumption on 
$\Gam{\setminus}(-\Lam)$. Hence $J^{*}(\del_{g})$ can be identified
with $e_{g}\in \kle(G)$, and we conclude that (N1) and (N2) are
satisfied. Moreover, since the $e_{g}$ are linearly independent and
$G$, being a compact infinite group, does not contain any isolated
point, (N3$'$) is fulfilled as well. It is left to apply
Lemma~\ref{2.3}, Remark~\Remark\  and Lemma~\ref{2.5}.
\eop
%\bigskip
%\addtocounter{theo}{1}
%
%\noindent
%{\em Example~\Examplethree}\/ \ 
%\label{Examplethree}%
%If $\Lam\subset \Z$ is a subset such that $\Z{\setminus}(-\Lam)$ is a
%Riesz set and if $\mu$ is a continuous measure on $\T$, then the
%convolution operator $T\dopu f\mapsto f\ast \mu$ on $C_{\Lam}(\T)$
%satisfies the \DE~{\rm(\ref{eqx.y})}. If $\mu$ is singular, then $T$
%is not weakly compact. In fact, as in Example~\Exampletwo\ 
%we have $\pi_{t}(q_{s})=0$  for all
%$t\in\T$ so that (\ref{eq2.3}) is fulfilled.
%If $\mu $ is continuous and
%singular, the proof that $T$ is not weakly compact
%is the same as in Example~\Exampletwo\ which
%deals with the special case $\Lam=\N$.
%
%Actually, the result extends to infinite compact abelian groups
%instead of the circle group.
\bigskip

We remark that some restriction on $\Lam$ is necessary in order to
ensure the \DE\ for -- say -- compact operators on $C_{\Lam}$,
because for the class of Sidon sets $\Lam$ the spaces $C_{\Lam}$ are
isomorphic to $\kle(\Lam)$, and thus there are one-dimensional
operators on $C_{\Lam}$ which fail the \DE\ \cite[Cor.~1]{Woj92}.

We finish this section with another negative result which is a
counterpart of the one just quoted.

\begin{prop}\label{3.8}
If $X^{*}$ has the Radon-Nikod\'ym property (in particular, if $X^{*}$
is separable), then there is a one-dimensional operator on $X$
failing the \DE.
\end{prop}

\Proof
Since $X^{*}$ has the Radon-Nikod\'ym property, $B_{X^{*}}$ contains
a weak$^{*}$ strongly exposed point $x_{0}^{*}$
\cite[Th.~5.12]{PheLNM}, that is, there is $x_{0}\in X$ such that
$\re x_{0}^{*}(x_{0}) = \|x_{0}^{*}\| = \|x_{0}\| = 1$ and
\beq\label{eq3.1}
\|x_{n}^{*}\| \le1 ,\ \re x_{n}^{*}(x_{0}) \to 1
\quad\Rightarrow\quad \|x_{n}^{*}-x_{0}^{*}\|\to 0.
\eeq
Define $T\in L(X)$ by $T(x)= x_{0}^{*}(x) x_{0}$ and assume that
$\|\Id-T\| = 1+\|T\| =2$. Then $\|x_{n}^{*}- T^{*}x_{n}^{*}\|\to 2$
for some sequence $(x_{n}^{*})\subset B_{X^{*}}$ and thus
$\|T^{*}x_{n}^{*}\|\to1$. Hence $|x_{n}^{*}(x_{0})|\to 1$  and with
no loss in generality $x_{n}^{*}(x_{0})\to\alp$ for some $|\alp|=1$.
This implies
$(\alp^{-1}x_{n}^{*})(x_{0})\to1$ and by (\ref{eq3.1})
$\|\alp^{-1}x_{n}^{*} - x_{0}^{*}\|\to0$. Finally we obtain the
contradiction
$$
2 = \lim_{n\to\iy} \|x_{n}^{*} - T^{*}x_{n}^{*}\| =
\|\alp x_{0}^{*} - T^{*}(\alp x_{0}^{*})\| =0. \eqno\Box
$$

\bigskip
In the setting of harmonic analysis this result tells us that for
many sets $\Lam\subset \Gam$, in particular for Shapiro sets
\cite{G-Riesz}, there are one-dimensional operator on $C(G)/C_{\Lam}$
failing the \DE, since for those sets $C(G)/C_{\Lam}$  has a
separable dual.

%
%
%  References
%
%
\typeout{References}
%\bibliography{az}
%\bibliographystyle{standard}

%
%
% Address
%
%
\small
\bigskip
\noindent
I.~Mathematisches Institut, Freie Universit\"at Berlin,
Arnimallee 2--6, \\
D-14\,195 Berlin, Germany; \ 
e-mail: 
werner@math.fu-berlin.de

\end{document}